\documentclass{article}

\usepackage{amssymb,dsfont,latexsym}

\newtheorem{thm}{Theorem}

\newtheorem{lem}[thm]{Lemma}

\newcommand\enu[1]{\smallskip\newline\makebox[5mm][l]{\rm(#1)}}

\newcommand\bp{\noindent{\it Proof.}\ }

\begin{document}
\title{The analogue of Choi matrices for a class of linear maps on Von Neumann algebras}

\author{Erling St{\o}rmer}

\date{12-18-2014}

\maketitle


 \section {Introduction}

 In the theory of positive linear maps of finite dimensional operator algebras Choi matrices  introduced by Choi \cite{C} play a central role, as they reduce the study of maps to that of concretely defined matrices.  In infinite dimensions a generalization of Choi matrices and their construction has been obtained by Holevo \cite{H} for maps of $B(H)$.  In the present paper we shall by a different approach extend the definition to maps of the form $C \to \sum_{i=1}^k  A_i CB_i$ on a Von Neumann factor $\mathcal{M}$ with a separating and cyclic vector, where $A_i, B_i, C \in \mathcal{M} $. We then show the analogue of the basic theorems for Choi matrices obtained in finite dimensions.
 
 The paper is divided into three sections.  The first, Section 2, treats factors with separating and cyclic vectors and is essentially a new and extended version of the operator part of \cite{BS}, but the proofs and emphasis are different.  In the second section, Section 3, we reformulate the basic part of the theory of Choi matrices so it fits into the setting of Section 2.  Then in Section 4 we define the analogue of Choi matrices for factors and prove theorems which extend the classical results for Choi matrices and completely positive maps
 
 \section {A spectral theorem}
 
 Throughout this section $\mathcal{M}$ is a factor acting on a Hilbert space $H$ and $x \in H$ is a unit vector which is separating and cyclic for $\mathcal{M}$.  We denote by $E$ the minimal projection in $B(H)$ onto the subspace spanned by $x$, and by $\omega$ the vector state $\omega(A) = (Ax,x)$. We then have the identity
 $$
 \omega(A)E = EAE   \     for A\in B(H).
 $$
 Let
 $$
 \mathcal{C} = \{ \sum_{i=1}^k A_iEB_i : A_i, B_i \in \mathcal{M}, k\in \mathbb{N}\}.
 $$
 Then $\mathcal{C}$ consists of finite rank operators in $B(H)$.

 \begin{lem}
$\mathcal{C}$ is a weakly dense *-algebra of $B(H)$.
\end{lem}
\bp Clearly $\mathcal{M}$ is a linear subspace of $B(H)$ and is self-adjoint since $(AEB)^* = B^*EA^* \in \mathcal{C} $ if $A,B \in \mathcal{M}$.  Furthermore, if $A,B,C,D \in \mathcal{M}$ then
$$
(AEB)(CED) = A(EBCE)D = \omega(BC)AED \in \mathcal{C},
$$
so by linearity $\mathcal{C}$ is closed under multiplication, hence is a *-subalgebra of $B(H)$.

To show density let first $y=Ax \in \mathcal{M}x$.  Since $E$ is of rank 1, so is $AEA^*$, since $x$ is separating for $\mathcal{M}$.  Thus $AEA^* = \lambda F$ for a rank 1 projection $F \in \mathcal{C}$ with $\lambda \in \mathbb{C}$. Then there exists $\mu \in \mathbb{C}$ such that 
$$
Fy = \lambda^{-1} AEA^*y = \lambda^{-1} A \mu x = \lambda^{-1} \mu y.
$$
hence y belongs to the range of $F$.

Let now $y$ be an arbitrary unit vector in $H$.  Since $x$ is cyclic for $\mathcal{M}$ there is a sequence $(x_i)$ in $\mathcal{M}x$ such that $x_i \to y$ in norm.  By the above paragraph for each $i$ there is a rank 1 projection $F_i \in \mathcal{C}$ such that $F_i   x_i  = x_i$.  Let $z\in H$ and $F$ as above. Then
$$
Fz = (z,y)y = \lim_i (z,x_i)x_i = \lim_i F_i z.
$$
Thus $F_i \to F$ strongly, so $F$ belongs to the weak closure of $\mathcal{C}$. Thus each rank 1 projection belongs to the weak closure of $\mathcal{C}$, hence by linearity so does each finite rank operator.  Thus $\mathcal{C}$ is weakly dense in $B(H)$. The proof is complete.
\medskip

Our next aim is to prove a spectral theorem for self-adjoint operators in $\mathcal{C}$.  For this we need two lemmas.

\begin{lem}
Let $P$ be a non-zero projection in $\mathcal{C}$.  Then there exists a rank 1 projection $F \in \mathcal{C}$ majorized by P.
\end{lem}
\bp  Since $x$ is cyclic for $\mathcal{M}$ there exists $A \in \mathcal{M}$ such that $PAx \neq 0$ or equivalently $PAE \neq 0$. Thus $PAEA^*P \neq 0$.  But this operator is in $\mathcal{C}$, has rank 1 and is positive, hence is a scalar multiple of a rank 1 projection $F \leq P$, proving the lemma.

\begin{lem}
Let $P \in \mathcal{C}$ be a rank 1 projection.  Then there exists $S \in \mathcal{C}$ such that $P = SES^*$.
\end{lem}
\bp  Let $y$ be a unit vector so that $Py = y$.  Since $P \in \mathcal{C}$ there are $A_i, B_i \in \mathcal{M}$ such that $P = \sum_i  A_i EB_i. $  Then $A_i EB_i y = \lambda_i A_i x$ for some $\lambda_i \in \mathbb{C}$ for all $i$. Let $S = \sum \lambda_i A_i \in \mathcal{M}$.  Then
$$
y = Py = \sum \lambda_i  A_i x = Sx.
$$
The operator $SES^*$ is positive of rank 1 in $\mathcal{C}$, so equals $\lambda Q$ for a rank 1 projection $Q \in \mathcal{M}$.  Note that $y \neq 0$ because if $Tr$ is the usual race on $B(H)$ then
$$
\lambda = Tr(SES^*) = Tr(ES^*S) = \omega(S^*S) \neq 0,
$$
since $S \in \mathcal{M}$ and $\omega$ is faithful on $\mathcal{M}$.  Let $z\in H$.  Then
$SES^*z = S \mu x$ with $\mu \in \mathbb{C}$.  Thus it follows, since $y = Sx$, that the range of $Q$ is the 1-dimensional subspace spanned by $y$, hence $Q= P$.  Hence if we scale $S$ we have $SES^* = P$.  The proof is complete.
\medskip

We can now prove the announced spectral theorem.

\begin{thm}\label{thm}
Let $T\in \mathcal{C}$ be self-adjoint.  Then each spectral projection of $T$ belongs to $\mathcal{C}$. Hence there exist $S_j \in \mathcal{M}$  and $c_j \in \mathbb{R}$ such that the operators $S_jES_j^*$ form an orthogonal family of rank 1 projections, and $T = \sum_j c_j S_j ES_j^*$.
\end{thm}
\bp Let $T = \sum_i A_i EB_i$.  Then
$$
T^2 = \sum_{ij} A_i EB_i A_j EB_j = \sum_{ij} \omega(B_iA_j)A_i EB_j.
$$
Iterating we get
$$
T^n = \sum_{ij} c_{ij}^{(n)} A_i EB_j,  \  c_{ij}^{(n)} \in \mathbb{C}.
$$
Hence each polynomial $p(T)$ in $T$ is of the form
$$
p(T) =  \sum_{ij}  d_{ij}^{(n)} A_i EB_j,  \  d_{ij}^{(n)} \in \mathbb{C},
$$
so that $p(T)$ belongs to the finite dimensional subspace $\mathcal{D}$ of $\mathcal{C}$ spanned by the operators $A_i EB_j$. Since $T$ is of finite rank each spectraøl projection of $T$ is a norm limit of polynomials in $T$, which by the above formula for $p(T)$ belongs to $\mathcal{D}$.  Since $\mathcal{D}$ is norm closed, each spectral projection of $T$ therefore belongs to $\mathcal{D}$, hence to $\mathcal{C}$.  Applying Lemmas 2 and 3 we obtain the conclusion of the theorem.  

\section{Choi matrices}

Let $\phi \colon M_n \to M_m$, where $M_n = M_{n} (\mathbb{C})$ denotes the complex $n \times n$ matrices, be a linear map.  Let $k \geq \max (n,m)$. We can then consider $M_n$ and $M_m$ as blocks in $M_k$ by choosing for $M_k$ matrix units $e_{ij}$ such that $e_{ij} \in M_n$ if and only if $i,j \leq n$, and similarly for $M_m$.  We can thus consider $\phi$ as a linear map $\bar \phi \colon M_k \to M_k$ by defining it by $\bar \phi(e_{ij}) = \phi(e_{ij}) $ if $i,j \leq n$ and 0 otherwise. We can thus restrict attention to linear maps $\phi \colon M_n \to M_n$.  Let $e_i,...,e_n$ be an orthonormal basis for $\mathbb{C}^n$ and $e_{ij}$ matrix units such that $e_{ij} e_k = \delta_{jk} e_i$.  Let $x = \sum_{i}  e_i  \otimes e_i  \in  \mathbb{C}^n  \otimes  \mathbb{C}^n$.  Then $x$ is separating and cyclic for $1 \otimes M_n \subset M_n \otimes M_n$, and 
$$
E = 1/n \sum_{ij} e_{ij} \otimes e_{ij}
$$
is the rank 1 projection onto $\mathbb{C}x$.  The Choi matrix for $\phi$ is defined to be
$$
C_{\phi} = \sum_{ij} e_{ij} \otimes \phi(e_{ij}) = n \iota_n \otimes \phi(E),
$$
where $\iota_n$ is the identity map on $M_n$. If $\phi$ has the form $\phi(C) = \sum_k A_kCB_k,$  with $A_k, B_k \in M_n$ then
$$
C_\phi = \sum_{ijk} e_{ij} \otimes A_k e_{ij} B_k = n \sum_k (1\otimes A_k)E(1\otimes B_k).
$$
The adjoint map for $\phi$ is defined by the formula
$$
Tr(\phi(A)B) = Tr(A\phi^*(B)).
$$
Thus with $\phi$ as above $\phi^*(C) = \sum B_k CA_k$, and 
$$
C_{\phi^*} = n \sum (1\otimes B_k)E(1\otimes A_k).
$$
By \cite{S} Lem. 4.1.10, $C_{\phi^*} = JC_{\phi}J $ for a conjugation $J$ of $\mathbb{C}^n \otimes \mathbb{C}^n$.  In particular $C_{\phi^*} \geq 0$ if and only if $C_{\phi} \geq 0$.  By the classical theorem of Choi \cite{C} $\phi$ is completely positive if and only if $C_{\phi} \geq 0.$  Hence $\phi$ is completely positive if and only if $C_{\phi^*} \geq = 0$, hence if and only if $\sum (1\otimes B_k) E(1\otimes A_k) \geq 0.$.

The main results on Choi matrices, which we shall generalize to factors in the next section, are the following two theorems, see \cite{S} Theorems 4.1.8 and 4.2.12. Recall that $\phi$ is positive, written $\phi \geq 0$ if $\phi(C) \geq 0$ whenever $C \geq 0$.

\begin{thm}\label{thm}
Let $\phi\colon M_n \to M_n$ be the map $\phi(C) = \sum A_i CB_i  \ with A_i,B_i \in M_n$.  Then the following conditions are equivalent:
\enu{i} $\phi$ is completely positive.
\enu{ii}  $\iota_n \otimes \phi \geq 0$.
\enu{iii} There are matrices $V_k \in M_n$ such that $\phi(C) = \sum_{k}  V_{k} ^* CV_k $.
\enu{iv} $C_{\phi^*} \geq 0.$
\enu{v}  $C_{\phi} \geq 0$.
\end{thm}

\begin{thm}\label{thm}
$\phi \geq 0$ if and only if for all positive $C,D\in M_n, Tr(C_{\phi} C \otimes D) \geq 0$.
\end{thm}

\section{Maps on factors}

We saw in the previous section that if $\phi(C) = \sum A_i CB_i$ with $A_i, B_i, C \in M_n$  the Choi matrix for $\phi$ is $n \sum  (1\otimes A_i )E(1\otimes B_i)$ and that of $\phi^*$ is $ n \sum (1\otimes B_i )E(1\otimes A_i)$.  In the present section we shall prove the analogous results of Theorems 5 and 6 for similar maps on factors.  The setting and notation will be the same as in section 2, so $\mathcal{M}$ is a factor acting on a Hilbert space $H$ with a separating and cyclic vector $X$, and $E$ is the projection onto $\mathbb{C}x$.  Let $\phi\colon \mathcal{M} \to \mathcal{M}$  be defined by
$$
\phi(C) = \sum_i A_i CB_i,  \  A_i, B_i \in \mathcal{M}.
$$
We denote by
$$
D_{\phi} = \sum B_i EA_i \in \mathcal{C},
$$
Where $\mathcal{C}$ is as in section 2.  Then $D_{\phi}$ corresponds to the Choi matrix $C_{\phi^*} $of section 3. The next result is the analogue of the Jamiolkowski isomorphism \cite{J}.

\begin{lem} \label {lem}
The map $\phi \to D_{\phi}$ with $\phi$ as above, is a linear isomorphism.
\end{lem}

\bp  The map is clearly linear.  To show injectivity assume $\phi = 0$.  Then for all $C \in \mathcal{M}$ and $D^{\prime}$ in the commutant $\mathcal{M}^{\prime}$ of $\mathcal{M}$ we have
$$
0 = Tr(\sum A_iCB_i D^{\prime}E) = Tr( CD^{\prime} \sum B_i EA_i).
$$
Since this holds for all $C\in \mathcal{M}, D^{\prime} \in \mathcal{M}^{\prime}$,   it holds for their linear combinations, hence for the Von Neumann algebra generated by $\mathcal{M}$ and $\mathcal{M}^{\prime}$.  Since $\mathcal{M}$ is a factor this algebra equals $B(H)$, hence by weak continuity of $\phi$ we have $\sum B_iEA_i =0 $.

Conversely, if $\sum B_i EA_i =0$ then for $C, D^{\prime}$ as above, we have 
\begin{eqnarray*}
0& =& Tr(D^{\prime} C \sum B_i EA_i) \\
&=& Tr(D^{\prime} \sum A_i CB_i E) \\
&=& (D^{\prime} \sum A_i CB_i x, x)  \\
&=& (\sum A_i CB_i x, D^{{\prime}*} x).
\end{eqnarray*}
Since $x$ is cyclic for $\mathcal{M}^{\prime}$ , $\sum A_i CB_i x= 0$.  Again, since $x$ separating for $\mathcal{M}$, $\sum A_i CB_i = 0$ for all $C \in \mathcal{M}$, i.e. $\phi = 0$.  It follows that $\phi = 0$ if and only if $D_{\phi} = 0$.  The proof is complete.

\begin{thm}\label{thm}
With $\phi$ as above, $D_{\phi} \geq 0$ if and only if there are $V_j \in \mathcal{M}$ such that $\phi(C) = \sum V_{j}^* C V_j$ for $C \in \mathcal{M}$.
\end{thm}
\bp Assume $D_{\phi} \geq 0$.  By Theorem 4 $D_{\phi} = \sum c_j E_j$, where $E_j = S_j ES_{j}^*$ are mutually orthogonal rank 1 projections in $\mathcal{C}$, and $c_j > 0$.  Let $ V_j  = c_{j}^{1/2} S_j$.  Then $D_{\phi} = \sum V_j E_j V_{j}^*$ , so by Lemma 7 $\phi(C) = \sum V_{j}^* C V_j$ for $C \in \mathcal{M}$.

Conversely, if $\phi(C)= \sum V_{j}^* CV_j$ then $D_{\phi} = \sum V_j EV_{j}^*$ is positive in $\mathcal{C}$. The proof is complete.
\medskip

The above theorem shows that if $D_{\phi} \geq 0$ then $\phi$ is completely positive. For the converse we need that $\mathcal{M}$ is injective.  The theorem will be stated as a direct generalization of Theorem 5.

\begin{thm}\label{thm}
Let $\mathcal{M}$ be an injective factor and $\phi \colon \mathcal{M} \to \mathcal{M}$ be of the form $\phi(C) = \sum A_i C B_i, A_i, B_i \in \mathcal{M}$. Then the following conditions are equivalent.
\enu{i} $\phi$ is completely positive.
\enu{ii} The map $C \to \sum A_i CB_i$ is positive on $B(H)$.
\enu{iii}  There are operators $V_j \in \mathcal{M}$ such that $\phi(C) = \sum V_{j}^* C V_j$ for all $C \in \mathcal{M}$.
\enu{iv} $D_{\phi} \geq 0$.
\enu{v}. $C_{\phi }\geq 0$.
\end{thm}
\bp We first show the equivalence of the first 4 conditions.  For the implications $(ii) \Rightarrow (iii) \Leftrightarrow (iv)$ and $(iii) \Rightarrow (i)$ we do not need to assume $\mathcal{M}$ is injective.  Assume (ii) holds.  Then for all $C\geq 0$ in $B(H)$ 
$$
0 \leq Tr(\sum A_i CB_i E) = Tr(C \sum B_i EA_i) ,
$$
Hence $D_{\phi} \geq 0,$ which by Theorem 8 is equivalent to (iii). Clearly $(iii) \Rightarrow (i).$ Thus it remains to show

$(i)\Rightarrow (ii)$. Assume $\phi$ is completely positive.  Since $\mathcal{M}$ is injective, so is $\mathcal{M}^{\prime}$, hence there exists an increasing sequence
$$
 \mathcal{N}_1 \subset \mathcal{N}_2 \subset ...
 $$
 of finite type I factors with union weakly dense in $\mathcal{M}^{\prime}$. Since $\mathcal{M}$ is a factor, the von Neumann algebra $\mathcal{M} \vee \mathcal {N}_n $ generated by $\mathcal{M}$ and $\mathcal{N}_n$ is isomorphic to $ \mathcal{M} \otimes \mathcal{N}_n$ via the map $CD \to C \otimes D$, where $C \in \mathcal{M}, D \in \mathcal{N}_n$.  Since $\phi$ is completely positive, $\phi \otimes \iota_n$ is positive on $\mathcal{M} \otimes \mathcal{N}_n$, hence the map $\phi(CD^{\prime}) = \phi(C)D^{\prime}$ defines a positive map on $\mathcal{M} \vee \mathcal{N}_n$ for all n.  Thus the extension $\bar\phi$ on $\mathcal{M} \vee \mathcal{M}^{\prime} = B(H)$ defined by $\bar\phi(CD^{\prime}) = \phi(C)D^{\prime}$ is positive, i.e. $\sum A_i CB_i \geq 0$ for all positive $C \in B(H)$, i.e. (ii) holds.  
 
 Finally we show $(iv) \Leftrightarrow (v)$. By the above if $D_\phi \geq 0$ then $(ii)$ holds, hence $ \sum A_i CB_i  \geq 0$ for all positive $C \in B(H)$, hence in particular for $C = E$, so $C_{\phi} = \sum A_i EB_i \geq 0$, and $(v)$ holds.
 
 Conversely if $C_\phi \geq 0$ we apply the equivalence $(ii) \Leftrightarrow (iv)$ to the map $\phi^*(C) = \sum B_i CA_i$. Then $D_{\phi^*} = C_{\phi} \geq 0$.  Thus by $(ii)$ applied to $\phi^*$, we see that $\phi^*$ is positive on $B(H)$, hence in particular $D_{\phi} = \sum B_i EA_i  \geq 0$, so $(iv)$ holds.
 The proof is complete.

 \medskip
 It should be remarked that since $\mathcal{C}$ is weakly dense in $B(H)$ by Lemma 1, condition $(ii)$ in Theorem 9 can be replaced by the weaker condition that the map $C \to \sum A_i CB_i  $ is positive on $\mathcal{C}.$

 We complete this section by showing the announced generalization of  Theorem 6.
 
 \begin{thm}\label{thm}
 Let $\mathcal{M}$ be a factor and $\phi(C) = \sum A_i CB_i$ as in Theorem 9.  Then $\phi \geq 0$ if and only if $Tr(D_{\phi} CD^{\prime}) \geq 0$ for all positive operators $C \in \mathcal{M}, D^{\prime} \in \mathcal{M}^{\prime}.$
 \end{thm}
 \bp Write $D^{\prime}$ as  $D^{\prime} = S^* S$  with $S \in \mathcal{M}^{\prime}$, and let $C$ be positive in $\mathcal{M}$. Then we have:
 \begin{eqnarray*}
Tr(D_{\phi}CD^{\prime} ) &=& Tr( \sum B_i EA_i CD^{\prime}) \\
&=& Tr(E \sum A_i CB_i D^{\prime})  \\
&=& (\sum A_i CB_i D^{\prime} x, x)  \\
&=& (\sum A_i C B_i Sx, Sx),
 \end{eqnarray*}
 where $x$ is the unit vector such that $Ex=x$.  Since $x$ is cyclic for $\mathcal{M}^{\prime}$ it follows that $Tr(D_\phi CD^{\prime}) \geq 0$  for all positive $C \in \mathcal{M}, D^{\prime} \in \mathcal{M}^{\prime}$  if and only if $\sum A_i CB_i \geq 0$ for all $C \in \mathcal{M}$, proving the theorem.

Department of Mathematics,  University of Oslo, 0316 Oslo, Norway.

e-mail  erlings@math.uio.no

\end{document}